\documentclass{amsart}
\usepackage{amsthm}
\usepackage{amssymb}
\usepackage{amsmath}
\usepackage{amsmath,amscd}
\usepackage{color}
\usepackage{hyperref}
\usepackage[all,arc]{xy}
\usepackage{graphicx}
\graphicspath{ {images/} }

\newtheorem{thm}{Theorem}[section]

\theoremstyle{definition}
\newtheorem{defn}[thm]{Definition}

\newtheorem{exmp}[thm]{Example}

\theoremstyle{remark}
\newtheorem{rem}[thm]{Remark}

\usepackage{tikz}
\usepgflibrary{shapes.geometric}
\usetikzlibrary{arrows}

\tikzstyle{V}=[draw, fill =black, circle, inner sep=0pt, minimum size=1.5pt]
\tikzstyle{wV}=[draw, fill =white, circle, inner sep=0pt, minimum size=4.5pt]
\tikzstyle{bV}=[draw, fill =black, circle, inner sep=0pt, minimum size=4.5pt]
\tikzstyle{over}=[draw=white,double=black,line width=2pt, double distance=.5pt]

\newcommand{\nc}{\newcommand}
\nc{\rnc}{\renewcommand}
\nc{\bb}[1]{{\mathbb #1}}
\nc{\bbA}{\bb{A}}\nc{\bbB}{\bb{B}}\nc{\bbC}{\bb{C}}\nc{\bbD}{\bb{D}}
\nc{\bbE}{\bb{E}}\nc{\bbF}{\bb{F}}\nc{\bbG}{\bb{G}}\nc{\bbH}{\bb{H}}
\nc{\bbI}{\bb{I}}\nc{\bbJ}{\bb{J}}\nc{\bbK}{\bb{K}}\nc{\bbL}{\bb{L}}
\nc{\bbM}{\bb{M}}\nc{\bbN}{\bb{N}}\nc{\bbO}{\bb{O}}\nc{\bbP}{\bb{P}}
\nc{\bbQ}{\bb{Q}}\nc{\bbR}{\bb{R}}\nc{\bbS}{\bb{S}}\nc{\bbT}{\bb{T}}
\nc{\bbU}{\bb{U}}\nc{\bbV}{\bb{V}}\nc{\bbW}{\bb{W}}\nc{\bbX}{\bb{X}}
\nc{\bbY}{\bb{Y}}\nc{\bbZ}{\bb{Z}}
\nc{\mbf}[1]{{\mathbf #1}}
\nc{\bfA}{\mbf{A}}\nc{\bfB}{\mbf{B}}\nc{\bfC}{\mbf{C}}\nc{\bfD}{\mbf{D}}
\nc{\bfE}{\mbf{E}}\nc{\bfF}{\mbf{F}}\nc{\bfG}{\mbf{G}}\nc{\bfH}{\mbf{H}}
\nc{\bfI}{\mbf{I}}\nc{\bfJ}{\mbf{J}}\nc{\bfK}{\mbf{K}}\nc{\bfL}{\mbf{L}}
\nc{\bfM}{\mbf{M}}\nc{\bfN}{\mbf{N}}\nc{\bfO}{\mbf{O}}\nc{\bfP}{\mbf{P}}
\nc{\bfQ}{\mbf{Q}}\nc{\bfR}{\mbf{R}}\nc{\bfS}{\mbf{S}}\nc{\bfT}{\mbf{T}}
\nc{\bfU}{\mbf{U}}\nc{\bfV}{\mbf{V}}\nc{\bfW}{\mbf{W}}\nc{\bfX}{\mbf{X}}
\nc{\bfY}{\mbf{Y}}\nc{\bfZ}{\mbf{Z}}
\nc{\bfa}{\mbf{a}}\nc{\bfb}{\mbf{b}}\nc{\bfc}{\mbf{c}}\nc{\bfd}{\mbf{d}}
\nc{\bfe}{\mbf{e}}\nc{\bff}{\mbf{f}}\nc{\bfg}{\mbf{g}}\nc{\bfh}{\mbf{h}}
\nc{\bfi}{\mbf{i}}\nc{\bfj}{\mbf{j}}\nc{\bfk}{\mbf{k}}\nc{\bfl}{\mbf{l}}
\nc{\bfm}{\mbf{m}}\nc{\bfn}{\mbf{n}}\nc{\bfo}{\mbf{o}}\nc{\bfp}{\mbf{p}}
\nc{\bfq}{\mbf{q}}\nc{\bfr}{\mbf{r}}\nc{\bfs}{\mbf{s}}\nc{\bft}{\mbf{t}}
\nc{\bfu}{\mbf{u}}\nc{\bfv}{\mbf{v}}\nc{\bfw}{\mbf{w}}\nc{\bfx}{\mbf{x}}
\nc{\bfy}{\mbf{y}}\nc{\bfz}{\mbf{z}}

\nc{\mcal}[1]{{\mathcal #1}}
\nc{\calA}{\mcal{A}}\nc{\calB}{\mcal{B}}\nc{\calC}{\mcal{C}}\nc{\calD}{\mcal{D}}
\nc{\calE}{\mcal{E}} \nc{\calF}{\mcal{F}}\nc{\calG}{\mcal{G}}\nc{\calH}{\mcal{H}}
\nc{\calI}{\mcal{I}}\nc{\calJ}{\mcal{J}}\nc{\calK}{\mcal{K}}\nc{\calL}{\mcal{L}}
\nc{\calM}{\mcal{M}}\nc{\calN}{\mcal{N}}\nc{\calO}{\mcal{O}}\nc{\calP}{\mcal{P}}
\nc{\calQ}{\mcal{Q}}\nc{\calR}{\mcal{R}}\nc{\calS}{\mcal{S}}\nc{\calT}{\mcal{T}}
\nc{\calU}{\mcal{U}}\nc{\calV}{\mcal{V}}\nc{\calW}{\mcal{W}}\nc{\calX}{\mcal{X}}
\nc{\calY}{\mcal{Y}}\nc{\calZ}{\mcal{Z}}
\nc{\fA}{\frak{A}}\nc{\fB}{\frak{B}}\nc{\fC}{\frak{C}} \nc{\fD}{\frak{D}}
\nc{\fE}{\frak{E}}\nc{\fF}{\frak{F}}\nc{\fG}{\frak{G}}\nc{\fH}{\frak{H}}
\nc{\fI}{\frak{I}}\nc{\fJ}{\frak{J}}\nc{\fK}{\frak{K}}\nc{\fL}{\frak{L}}
\nc{\fM}{\frak{M}}\nc{\fN}{\frak{N}}\nc{\fO}{\frak{O}}\nc{\fP}{\frak{P}}
\nc{\fQ}{\frak{Q}}\nc{\fR}{\frak{R}}\nc{\fS}{\frak{S}}\nc{\fT}{\frak{T}}
\nc{\fU}{\frak{U}}\nc{\fV}{\frak{V}}\nc{\fW}{\frak{W}}\nc{\fX}{\frak{X}}
\nc{\fY}{\frak{Y}}\nc{\fZ}{\frak{Z}}
\nc{\fa}{\frak{a}}\nc{\fb}{\frak{b}}\nc{\fc}{\frak{c}} \nc{\fd}{\frak{d}}
\nc{\fe}{\frak{e}}\nc{\fFf}{\frak{f}}\nc{\fg}{\frak{g}}\nc{\fh}{\frak{h}}
\nc{\fri}{\frak{i}}\nc{\fj}{\frak{j}}\nc{\fk}{\frak{k}}\nc{\fl}{\frak{l}}
\nc{\fm}{\frak{m}}\nc{\fn}{\frak{n}}\nc{\fo}{\frak{o}}\nc{\fp}{\frak{p}}
\nc{\fq}{\frak{q}}\nc{\fr}{\frak{r}}\nc{\fs}{\frak{s}}\nc{\ft}{\frak{t}}
\nc{\fu}{\frak{u}}\nc{\fv}{\frak{v}}\nc{\fw}{\frak{w}}\nc{\fx}{\frak{x}}
\nc{\fy}{\frak{y}}\nc{\fz}{\frak{z}}

\nc{\al}{\alpha}
\nc{\ep}{\epsilon}
\nc{\la}{\lambda}
\nc{\be}{\beta}
\nc{\de}{\delta}
\nc{\na}{\nabla}
\nc{\nap}{{\na_+}}
\nc{\nam}{{\na_-}}
\nc{\lefta}{\leftarrow}
\nc{\wt}{\widetilde}

\DeclareMathOperator{\stab}{stab}

\DeclareMathOperator{\Attr}{Attr}
\DeclareMathOperator{\FAttr}{FAttr}
\DeclareMathOperator{\supp}{supp}

\DeclareMathOperator{\pt}{pt}

\DeclareMathOperator{\rank}{rank}

\DeclareMathOperator{\Lie}{Lie}

\DeclareMathOperator{\opp}{opp}

\DeclareMathOperator{\gr}{gr}
\DeclareMathOperator{\Ind}{Ind}

\newcommand{\hatQ}{\hat{Q}}
\newcommand{\unit}{\mathbf{1}}  
\DeclareMathOperator{\ev}{ev}

\newcommand{\taup}{{\overset{+}\tau}}
\newcommand{\taum}{{\overset{-}\tau}}
\newcommand{\taupm}{{\overset{\pm}\tau}}

\nc{\stp}[2]{\stab^{+, {#1}}_{#2}}   
\nc{\stm}[2]{\stab^{-,{#1}}_{#2}} 
\nc{\stabp}{\stab^+} 
\nc{\stabm}{\stab^-} 

\nc{\fstab}{\fs\ft\fa\fb}
\DeclareMathOperator{\ob}{{b}}
\DeclareMathOperator{\oR}{R}

\DeclareMathOperator{\aff}{aff}

\DeclareMathOperator{\Coh}{Coh}
\DeclareMathOperator{\Mod}{Mod}
\DeclareMathOperator{\Frac}{Frac}
\DeclareMathOperator{\Rep}{Rep}
\DeclareMathOperator{\Char}{Char}

\newcommand{\Gm}{\mathbb{G}_m}
\newcommand{\La}{\Lambda}
\DeclareMathOperator{\Hom}{Hom}

\title{Stable bases of the Springer resolution and representation theory}
\date{\today}

\author[C.~Su]{Changjian~Su}
\address{Department of mathematics, University of Toronto, Toronto, ON, Canada}
\email{changjiansu@gmail.com}

\author[C.~Zhong]{Changlong~Zhong}
\address{State University of New York at Albany, 1400 Washington Ave, ES 110, Albany, NY, 12222}
\email{czhong@albany.edu}

\begin{document}

\maketitle
\begin{abstract} 
In this note, we collect basic facts about Maulik and Okounkov's stable bases for the Springer resolution, focusing on their relations to representations of Lie algebras over $\bbC$ and algebraically closed positive characteristic fields, and of the Langlands dual group over non-Archimedean local fields. 
\end{abstract}


\section{Introduction}
The Schubert calculus studies the cohomology and K-theory of the (partial) flag varieties $\calB$. In this note, we study certain basis elements of Maulik and Okounkov, called the \textit{stable basis}, in the equivariant cohomology and K-theory of the cotangent bundle of the flag varieties $T^*\calB$ (the Springer resolution). Pulling back to the zero section $\calB$, the stable basis becomes some natural basis in the equivariant cohomology and K-theory of the flag varieties.

Maulik and Okounkov introduced stable basis in their study of quantum cohomology of  quiver varieties \cite{MO}. Later, Okounkov and his collaborators introduced  K-theory and elliptic cohomology versions of these bases; see \cite{O, OS16, AO}. They turn out to be very useful both in enumerative geometry and geometric representation theory; see the surveys by Okounkov \cite{O15, O18}. 

The Maulik--Okounkov stable basis are constructed for a class of varieties called symplectic resolutions \cite{Kal09}. We will focus on the Springer resolution, which is one of classical examples of symplectic resolutions. Both the cohomological and K-theoretic stable basis for the Springer resolution are the standard objects for certain  category.

In the cohomology case, the stable basis are just some rational combination of the conormal bundles of the Schubert cells, which was shown  to coincide with  characteristic cycles of certain $\calD$-modules on the flag variety. Via the localization theorem of Beilinson and Bernstein \cite{BB}, these $\calD$-modules correspond to the Verma modules for the Lie algebra. Furthermore, the action of the convolution algebra of the Steinberg variety \cite{CG}, which is isomorphic to the graded affine Hecke algebra \cite{L88}, is computed for the stable basis \cite{Su17}. From this, the first author deduced the localization formula for the stable basis, which are direct generalizaiton of the well-known AJS/Billey formula \cite{B99} for the localization of  Schubert classes in the equivariant cohomology of flag varieties. The restriction formulae also play a crucial role in determining the quantum connection of the cotangent bundle of partial flag varieties \cite{Su17}.

On the other hand, the graded affine Hecke algebra also appears in the work of Aluffi and Mihalcea \cite{AM16} about the Chern-Schwartz-MacPherson (CSM) classes \cite{M, S65a, S65b} for  Schubert cells. Comparing the actions, we identity the pullbacks of the stable basis with the CSM classes \cite{AMSS17} (see also \cite{RV15}). The effectivity of the characteristic cycle of $\calD$-modules enables the authors in \cite{AMSS17} to prove the non-equivariant positivity conjecture of Aluffi and Mihalcea \cite{AM16} for CSM classes of the Schubert cells.

The K-theoretic case is even more interesting. By the famous theorem of Kazhdan, Lusztig and Ginzburg ( see \cite{KL, G, CG}), the affine Hecke algebra is isomorphic to the convolution algebra of the equivariant K-theory of the Steinberg variety. Hence, the affine Hecke algebra acts on the equivariant K-theory of $T^*\calB$. This action on the stable basis is computed in \cite{SZZ17}. It roughly says that the stable basis is the standard basis for the regular representation of the finite Hecke algebra. With this action, we can compute the localization of the stable basis via the root polynomial method. 

The K-theretic generalization of the CSM classes are motivic Chern classes defined by Brasselet, Sch{\"u}rmann and Yokura \cite{BSY}. Pulling back the K-theory stable basis to the zero section, we get the motivic Chern classes of the Schubert cells; see \cite{AMSS19, FRW}. By definition, the motivic Chern classes enjoy good functorial properties. These facts enable the authors in \cite{AMSS19} to use Schubert calculus to prove some conjectures of Bump, Nakasuji and Naruse (\cite{BN17, N14}) about unramified principle series of the $p$-adic Langland dual group. 

Recall the affine Hecke algebra can be realized either as the equivariant K-group of the Steinberg variety, or as the double Iwahori-invariant functions on the $p$-adic Langland dual group; see \cite[Introduction]{CG} and \cite{B16}. The bridge connecting the Schubert calculus over complex numbers and the representation theory of the $p$-adic Langland dual group is the shadow of these two geometric realizations. To be more precise, the (specialized) equivariant K-theory of the flag variety and the Iwahori-invariants of an unramified principle series of the $p$-adic group are two geometric realizations of the regular representation of the finite Hecke algebra. Under this isomorhpism, the stable basis (or the motivic Chern classes of the Schubert cells) gets identified with the standard basis on the $p$-adic side, while the fixed point basis maps to the Casselman basis \cite{C80}. With this isomorphism, we can give an equivariant K-theory interpretation of the Macdonald's formula for the spherical function (\cite{C80}) and the Casselman--Shalika formula (\cite{CS}) for the spherical Whittaker function (see \cite{SZZ17}).

The stable basis is also related to representations of Lie algebras in positive characteristic fields, via the localization theorem of Bezrukavnikov, Mirkovi\'{c} and Rumynin \cite{BMR1, BMR2}, which generalizes the famous localization theorem of Beilinson and Bernstein \cite{BB} over the complex numbers. This is achieved as follows. The definition of the K-theory stable basis depends on a choice of the alcove. Changing alcoves defines the so-called wall R-matrices (\cite{OS16}). We first recall these wall R-matrices for the Springer resolution, which are computed in \cite{SZZ19}. The formulae can be nicely packed using the Hecke algebra actions. It turns out the wall R-matrices coincide with the monodromy matrices of the quantum connection. Bezrukavnikov and Okounkov conjecture that for general symplectic resolutions $X$, the monodromy representation is isomorphic to the representation coming from derived equivalences (\cite{O15, BO}). The Springer resolution case is established by Braverman, Maulik and Okounkov in \cite{BMO}. Thus, the wall R-matrices is also related to the derived equivalences. Finally, we give a categorification of the stable basis via the affine braid group action, constructed by Bezrukavnikov and Riche (see \cite{BR, R08}), on the derived category of coherent sheaves on the Springer resolution. Moreover, the categorified stable basis is identified with the Verma modules for the Lie algebras over positive characteristic fields under the localization equivalence. 

\subsection*{Acknowledgments}
We thank G. Zhao for the collaboration of \cite{SZZ17, SZZ19}. The first-named author also thanks P. Aluffi, L. Mihalcea and J. Sch{\"u}rmann for the collaboration of \cite{AMSS17, AMSS19}. We thank the organizers for invitation to speak at the `International Festival in Schubert Calculus', Guangzhou, China, 2017.

\section{Cohomological Stable basis}\label{sec:cohcase}
Let us first recall the definition of the stable basis of Maulik and Okounkov \cite{MO}. 

\subsection{Notations}
Let $G$ be a connected, semisimple, complex linear algebraic group with a Borel subgroup $B$ and a maximal torus $A$. Let $B^-$ be the opposite Borel subgroup. Let $\La$ (resp. $\La^\vee$) be the group of characters (resp. cocharacters) of $A$. Let $R^+$ denote the roots in $B$ and $Q$ be the root lattice. Let $\fg$ (resp. $\fh$) be the Lie algebra of $G$ (resp. $A$). Let $\rho$ be the half-sum of the positive roots. Let $\fC_\pm\subset \fh$ be the dominant/anti-dominant Weyl chamber in $\fh$. Let $\calB=G/B$ be the variety of all Borel subgroups of $G$. Denote $T^*\calB$ (resp. $T\calB$) the cotangent (tangent) bundle. The cotangent bundle $T^*\calB$ is a resolution of the nilpotent cone $\calN\subset \fg^*$, which is the so-called Springer resolution. The group $\bbC^*$ acts on $T^*\calB$ by  $z\cdot (B', x)=(B', z^{-2}x)$ for any $z\in \bbC^*, (B', x)\in T^*\calB$. Let $-\hbar$ be the $\Lie \bbC^*$-weight of the cotangent fiber. Let $T=A\times\bbC^*$. Then $H_T^*(pt)=\bbC[\Lie A][\hbar]$. Let $W$ be the Weyl group, and let $\leq$ denote the Bruhat order on $W$. The $A$-fixed points in $T^*\calB$ are indexed by the Weyl group $W$. For each $w\in W$, the corresponding fixed point is $wB\in \calB\subset T^*\calB$. Let $\iota_w:wB \hookrightarrow T^*\calB$ denote the embedding. For any $\gamma\in H_T^*(T^*\calB)$, let $\gamma|_w=\iota_{w}^*\gamma\in H_T^*(wB)=H_T^*(pt)$.

\subsection{Definition}
The definition of the stable basis depends on a choice a Weyl chamber $\fC$ in $\fh$. Pick any one parameter $\sigma:\bbC^*\rightarrow A$ lying in the chamber $\fC$. Recall the attracting set of the torus fixed point $wB\in T^*\calB$
\[
\Attr_\fC(w)=\{x\in T^*\calB|\underset{z\to 0}\lim \sigma(z)\cdot x=wB\}.
\]
A partial order on $W$ is defined as
\[
v\preceq_\fC w\quad \text{if }\quad \overline{\Attr_\fC(w)}\cap vB\neq \emptyset.
\]
The order determined by the dominant chamber (resp. anti-dominant chambert) is the usual Bruhat order (resp., the opposite Bruhat order). The full attracting set is 
\[
\FAttr_\fC(w)=\cup_{v\preceq_\fC w}\Attr_\fC(v). 
\]

For each $w\in W$, let $\epsilon_w=e^A(T_{wB}^*\mathcal{B})$, the $A$-equivariant Euler class of the $A$-vector space $T_{wB}^*\mathcal{B}$. I.e., $\epsilon_w$ is the product of $A$-weights in the vector space $T_{wB}^*\mathcal{B}$. Let $N_w=T_{wB}(T^*\mathcal{B})$. The chamber $\mathfrak{C}$ gives a decomposition of the tangent bundle
\[
N_w=N_{w,+}\oplus N_{w,-}
\]
into $A$-weights which are positive and negative on $\mathfrak{C}$ respectively. Since $T=A\times \mathbb{C}^*$, $H_T^*(pt)=H_A^*(pt)[\hbar]$. The sign in $\pm e(N_{w,-})\in H_T^*(pt)$ is determined by the condition
\[
\pm e(N_{w,-})|_{\hbar=0}=\epsilon_w.
\]
The cohomological stable basis is defined as follows.
\begin{defn}\cite[Theorem 3.3.4]{MO}\label{stable theorem B}
There exists a unique map of $H_T^*(pt)$-modules
\[
\stab_{\mathfrak{C}}:H_T^*((T^*\mathcal{B})^A)\rightarrow H_T^*(T^*\mathcal{B})
\]
satisfying the following properties. For any $w\in W$, denote $\stab_{\mathfrak{C}}(w)=\stab_{\mathfrak{C}}(1_w)$. Then
\begin{enumerate}
\item (Support) $\supp \stab_{\mathfrak{C}}(w)\subset \FAttr_\fC(w)$,
\item (Normalization) $\stab_{\mathfrak{C}}(w)|_w=\pm e(N_{-,w})$, with sign according to $\epsilon_w$,
\item (Degree) $\stab_{\mathfrak{C}}(w)|_v$ is divisible by $\hbar$, for any $v\prec_{\mathfrak{C}} w$.
\end{enumerate}
\end{defn}
\begin{rem}
\begin{enumerate}
\item
The map is defined by a Lagrangian correspondence between $(T^*\mathcal{B})^A\times T^*\mathcal{B}$, hence mapping middle degree to middle degree. Therefore, the last condition is equivalent to 
\[\deg_A \stab_{\mathfrak{C}}(w)|_v<\deg_A\stab_{\mathfrak{C}}(v)|_v.\] 
Here, for any $f\in H_T(pt)=\bbC[\Lie A][\hbar]$, $\deg_A f$ denotes the degree of the polynomial $f$ in the $A$-variables.
\item 
By the first and  second conditions above, $\{\stab_{\mathfrak{C}}(w)| w\in W\}$ are a basis for the localized equivariant cohomology $H_T^*(T^*\calB)_{\textit{loc}}:=H_T^*(T^*\calB)\otimes_{H_T(pt)}\Frac H_T(pt)$. It is the so-called \textbf{stable basis}.
\item
The stable bases for opposite chambers are dual bases; see \cite[Theorem 4.4.1]{MO}
. I.e., 
\[\langle\stab_{\mathfrak{C}}(v), \stab_{\mathfrak{-C}}(w)\rangle=(-1)^{\dim G/B}\delta_{v,w},\]
where the pairing $\langle - , -\rangle$ is the Poincar\'e pairing defined via localization.
\end{enumerate}
\end{rem}
Let $\stab_+(w)=\stab_{\fC_+}(w)$ and $\stab_-(w)=\stab_{\fC_-}(w)$. 

\subsection{The graded affine Hecke algebra action}\label{sec:graded}
The graded affine Hecke algebra $\calH_\hbar$ is generated by the elements $x_\lambda$ for $\lambda\in \fh^*$, $w\in W$ and a central element $\hbar$ such that 
\begin{enumerate}
\item 
$x_\lambda$ depends linearly on $\lambda\in \fh^*$;
\item 
$x_\lambda x_\mu=x_\mu x_\lambda$ for any $\lambda, \mu\in \fh^*$;
\item 
The $w$'s form the Weyl group inside $\calH_\hbar$;
\item 
For any simple root $\alpha_i$ and $\lambda\in \fh^*$, we have
\[s_ix_\lambda-x_{s_i\lambda}s_i=\hbar (\lambda, \alpha_i^\vee).\]
\end{enumerate}
According to \cite{L88}, there is a natural algebra isomorphism 
\[\calH_\hbar\simeq H_*^{G\times \bbC^*}(Z),\]
where $Z=T^*\calB\times_\calN T^*\calB$ is the Steinberg variety, and the right hand side is endowed with the convolution algebra structure. The isomorphism is constructed as follows. For any $\lambda\in \La$, let $\calL_\lambda=G\times_B\bbC_\lambda$ be the line bundle on $\calB$. Pulling back to $T^*\calB$, we get a line bundle on $T^*\calB$, which is still denoted by $\calL_\lambda$. The the above isomorphism sends $x_\lambda$ to the push-forward of the first Chern class $c_1(\calL_\lambda)$ under the diagonal embedding $T^*\calB\hookrightarrow Z$. The image of a simple reflection $s_\alpha$ in the Weyl group is constructed as follows. Let $P_\alpha$ denote the corresponding minimal parabolic subgroup containing $B$. Let $\calP_\alpha=G/P_\alpha$ and $Y_\al=\calB\times_{\calP_\al}\calB\subset \calB\times \calB$. The conormal bundle $T^*_{Y_\al}(\calB\times \calB)$ is denoted by $T^*_{Y_\al}$, which is a smooth closed $G\times\bbC^*$-invariant subvariety of $Z$. Then $s_\alpha-1$ is sent to the cohomology class $[T^*_{Y_\al}]\in H_*^{G\times \bbC^*}(Z)$.

Via this isomorphism, $\calH_\hbar$ acts on $H_T^*(T^*\calB)$ by convolution (\cite[Chapter 2]{CG}). Let $\pi$ denote this action. The action on the stable basis is given by
\begin{thm}\cite[Lemma 3.2]{Su17}\label{thm:cohaction}
For any $w\in W$ and simple root $\alpha$,
\[\pi(s_\alpha)(\stab_\pm(w))=-\stab_\pm(ws_\alpha).\]
\end{thm}

\subsection{The restriction formula}
One immediate corollary of the above theorem is the restriction formula for the stable basis; see \cite[Theorem 1.1]{Su17}.
\begin{thm}\label{restriction -B}
Let $y=\sigma_1\sigma_2\cdots\sigma_l$ be a reduced expression for $y\in W$. Then
\begin{equation}\label{formula restriction -B}
\stab_-(w)|_y=(-1)^{l(y)}\prod\limits_{\alpha\in R^+\setminus R(y)}(\alpha-\hbar)\sum\limits_{\substack{1\leq i_1<i_2<\dots<i_k\leq l\\w=\sigma_{i_1}\sigma_{i_2}\dots\sigma_{i_k}}}\hbar^{l-k}\prod\limits_{j=1}^k \beta_{i_j},
\end{equation}
where $\sigma_i$ is the simple reflection associated to a simple root $\alpha_i$, $\beta_i=\sigma_1\cdots\sigma_{i-1}\alpha_i$, and $R(y)=\{\beta_i|1\leq i\leq l\}$. Furthermore, the sum in Equation (\ref{formula restriction -B}) does not depend on the reduced expression for $y$.
\end{thm}
From this formula, we can recover the AJS/Billey formula (\cite{B99})
\[[\overline{B^-wB/B}]|_y=\sum\limits_{\substack{1\leq i_1<i_2<\dots<i_k\leq l\\ 
w=\sigma_{i_1}\sigma_{i_2}\dots\sigma_{i_k} \text{reduced}}}\beta_{i_1}\cdots\beta_{i_k},\]
via the following limit formula
\[[\overline{B^-wB/B}]|_y=(-1)^{\ell(w)}\lim_{\hbar\rightarrow \infty}\frac{\stab_-(w)|_y}{(-\hbar)^{\dim\overline{B^-wB/B}}}.\]
See \cite[Theorem 4.8]{Su17} for more details.

\subsection{Relation with Chern-Schwartz-MacPherson (CSM) classes}
Let us first recall the definition of Chern-Schwartz-MacPherson classes. According to a conjecture attributed to Deligne and Grothendieck, there is a unique natural transformation $c_*: \calF \to H_*$ from the functor $\calF$ of constructible functions on a complex algebraic variety
$X$ to homology of $X$, where all morphisms are proper, such that if $X$ is smooth
then $c_*(1_X)=c(TX)\cap [X]$. This conjecture was proved by MacPherson \cite{M}. The class $c_*(1_X)$ for possibly singular $X$ was shown to coincide with a class defined
earlier by M.-H.~Schwartz \cite{S65a, S65b}. For any
constructible subset $W\subset X$, we call the class $c_{SM}(W):=c_*(1_W)\in H_*(X)$ the \textit{Chern-Schwartz-MacPherson} (CSM) class of $W$ in $X$. The theory of CSM classes was later extended to the equivariant setting by Ohmoto \cite{O06}. 

Let $X(w)^\circ=BwB/B\subset X$ and $Y(w)^\circ=B^-wB/B$ be the Schubert cells in $\calB$. Let $X(w)=\overline{X(w)^\circ}$ be the Schubert variety. Similar formula in Theorem \ref{thm:cohaction} was also obtained by Aluffi and Mihalcea for the CSM classes for the Schubert cells $c_{SM}(X(w)^\circ)$; see \cite{AM16}. Comparing these formulae, it is easy to get the following relation between the stable basis and the CSM classes; see \cite[Corollary 1.2]{AMSS17} and \cite{RV15}.
\begin{thm}\label{thm:cohpullback}
Let $i:\calB\hookrightarrow T^*\calB$ be the inclusion. For any $w\in W$, 
\[\iota^*(\stab_+(w))|_{\hbar=1}=(-1)^{\dim G/B}c_{SM}(X(w)^\circ)\in H_A^*(\calB),\]
and 
\[\iota^*(\stab_-(w))|_{\hbar=1}=(-1)^{\dim G/B}c_{SM}(Y(w)^\circ)\in H_A^*(\calB).\]
\end{thm}

The equivariant cohomology of the flag variety $\calB$ has a natural basis, namely, the Schubert basis $\{[X(w)]|w\in W\}$. Thus we can expand the CSM classes in terms of this basis
\[c_{SM}(X(w)^\circ)=\sum_{u}c(w,u)[X(u)]\in H_A^*(\calB),\]
where $c(w,u)\in H_A^*(pt)$. It is conjectured by Aluffi and Mihalcea that (\cite{AM16})
\[c(w,u)\in \bbZ_{\geq 0}[\alpha|\alpha>0].\]
The non-equivariant case of this conjecture is proved in \cite{AMSS17}, in which the relation between the stable basis and the characteristic cycles of homolomic $\calD_{\calB}$-modules played an important role.

\subsection{Characteristic cycles of $\calD$-modules}
For $w \in W$, let $M_w$ be the Verma module of highest weight $-w\rho - \rho$, a module over the universal enveloping algebra $U(\fg)$. Let $\mathfrak{M}_w$ denote the holonomic $\calD_{\calB}$-module 
\[ \mathfrak{M}_w = \calD_\calB \otimes_{U(\fg)} M_w \/. \] 
The image of this regular holonomic $\calD_\calB$-module under the Riemann--Hilbert correspondence is the constructible complex $\bbC_{X(w)^\circ}[\ell(w)]$; see \cite{BB, BK81}. The characteristic cycle $\Char(\fM_w)$ of $\fM_w$ is a linear combination of the conormal bundles of the Schubert cells. The relation with the stable basis is given by the following formula. It is claimed in \cite[Remark 3.5.3]{MO}, and later proved in \cite[Lemma 6.5]{AMSS17}.
\begin{thm}\label{thm:char}
For any $w\in W$, 
\[\stab_+(w)=(-1)^{\dim G/B-\ell(w)}[\Char(\fM_w)]\in H_T^*(T^*\calB).\]
\end{thm}

\section{K-theoretic Stable bases and the affine Hecke algebra action}\label{sec:hecke}
We first recall the definition of K-theoretic stable bases of Maulik and Okounkov. We follow the notations of \cite{SZZ17}. 

\subsection{Notations continued} We will follow the notations in Section \ref{sec:cohcase}. Let us introduce more notations.  Let $H_{\al^\vee,n}=\{\la\in \fh^*_\bbR|(\la, \al^\vee)=n\}$ be the hyperplanes determined by the coroot $\al^\vee$ and integers $n$. The union of the hyperplanes is a closed subset of $\fh^*_\bbR$, whose complement has  connected components  called alcoves. The fundamental alcove is $\na_+=\{\la\in \fh^*_\bbR|0<\langle \la, \al^\vee\rangle <1, \text{for any positive coroot }\al^\vee\}$. Denote $\na_-=-\na_+$. 


In the remaining parts of this note, we will consider the equivariant K-theory, a good introduction of which can be found in \cite{CG}. If a group $H$ acts on an algebraic variety $X$, then the $H$-equivariant K-theory of $X$, which is denoted by $K_H(X)$, is defined to be the Grothendieck group of the $H$-equivariant coherent sheaves on $X$. I.e., $K_H(X):=K^0(\Coh_H(X))$. By definition, $K_H(X)$ is a module over $K_H(pt)=K^0(\Rep(H))$.

Recall the group $\bbC^*$ acts on $T^*\calB$ by  $z\cdot (B', x)=(B', z^{-2}x)$ for any $z\in \bbC^*, (B', x)\in T^*\calB$. Let $q^{-1}$ be the $\bbC^*$-character of of the cotangent fiber under this action, i.e., $q=e^\hbar$. Therefore, the $\bbC^*$-equivariant K-group of a point is $K_{\bbC^*}=K^0(\Rep(\bbC^*))=\bbZ[q^{ 1/2}, q^{-1/2}]$. Recall $T=A\times\bbC^*$. Denote $R=\bbZ[q^{1/2}, q^{-1/2}]$, and $S=K_T(\pt)=\bbZ[q^{1/2}, q^{-1/2}][\La]$. We consider the $T$-equivariant K-theory of $T^*\calB$. Recall $\iota_w$ denotes the inclusion of the fixed point $wB$ into $T^*\calB$. For any $\calF\in K_T(T^*\calB)$, let $\calF|_w$ denote $\iota_w^*\calF\in K_T(pt)$. By localization theorem \cite[Chapter 5]{CG}, the localized K-group $K_T(T^*\calB)_{\textit{loc}}:=K_T(T^*\calB)\otimes_{K_T(pt)}\Frac K_T(pt)$ has a fixed point basis $\{\iota_{w*}1 | w\in W\}$.

The non-degenerate pairing on $K_T(T^*\calB)$ can be defined using localization as follows
\[
\langle\calF, \calG\rangle=\sum_{w\in W}\frac{\calF|_w\calG|_w}{\bigwedge^\bullet T_w(T^*\calB)}=\sum_{w\in W}\frac{\calF|_w\calG|_w}{\prod_{\al>0}(1-e^{w\al})(1-qe^{-w\al})}, \quad \calF, \calG\in K_T(T^*\calB).
\]
Here for each $T$-space $V$, $\bigwedge ^\bullet V=\sum_{k\ge 0}(-1)^k\bigwedge^kV^\vee=\prod (1-e^{-\al})\in K_T(\pt)$, where the product runs through all $T$-weights in $V$, counted with multiplicities. 

\subsection{Definition of the stable basis} 
A polarization $T^{1/2}$ is a solution of the following equation 
\[
T^{1/2}+q^{-1}(T^{1/2})^\vee=T(T^*\calB) \quad  \text{in }K_T(T^*\calB). 
\]
Denote $T^{1/2}_{\opp}=q^{-1}(T^{1/2})^\vee$. We will frequently focus on the  two mutually opposite polarizations:  $T\calB$ and $T^*\calB$. Recall $N_w=T_{wB}(T^*\calB)$ is the normal bundle at $w$, and the chamber $\fC$ determines a decomposition $N_w=N_{w, +}\oplus N_{w,-}$ according to the sign of $A$-weights with respect to the chamber $\fC$. Let $N_w^{1/2}=N_w\cap (T^{1/2}|_w)$, and similarly define $N_{w, +}^{1/2}$ and $N_{w, -}^{1/2}$. It follows that the square root $(\frac{\det N_{w,-}}{\det N_{w}^{1/2}})^{1/2}$ exists in $K_T(T^*\calB)$. 

For any Laurent polynomial $f=\sum_{\mu\in \La}f_\mu e^\mu\in K_T(\pt)$, with $f_\mu\in \bbZ[q^{1/2}, q^{-1/2}], e^\mu\in K_A(\pt)$, define the $A$-degree of $f$ to be the Newton polygon
\[
\deg_Af=\text{ Convex hull }(\{\mu|f_\mu\neq 0\})\subset \La\otimes \bbR. 
\]

We can now recall the definition of the K-theory stable basis. 
\begin{defn}\cite{O} For any chamber $\fC$, polarization $T^{1/2}$, and alcove $\na$, there is a unique map  of $K_T(\pt)$-modules (called the stable envelope):
	\[
	\stab_{\fC,T^{1/2},\na}: K_{T}((T^*\calB)^{A})\to K_T(T^*\calB),
	\]
	satisfying the following properties. Denote $\stab^{\fC,T^{1/2}, \na}_w=\stab_{\fC,T^{1/2},\na}(1_w)$. Then
	\begin{enumerate}
		\item (Support) $\supp( \stab^{\fC, T^{1/2},\na}_w)\subset \FAttr_{ \fC}(w)$.
		\item (Normalization) $\stab^{\fC,T^{1/2}, \na}_w|_w=(-1)^{\rank N^{1/2}_{w, +}}\left(\frac{\det N_{w,-}}{\det N_w^{1/2}}\right)^{1/2} \calO_{\Attr _\fC(w)}|_w.$
		\item (Degree) $\deg _A(\stab^{\fC, T^{1/2},\na}_w|_v)\subset \deg_A(\stab^{\fC,T^{1/2}, \na}_v|_v)+v\lambda- w\lambda$ for any $v\prec_\fC w, \lambda\in \na$. 
	\end{enumerate}
\end{defn}
Note that the the degree condition depends on the alcove $\na$ only, not on a particular $\la\in \na$. On the other hand, the normalization does not depend on the alcove. 

Denote 
\[\stabp_{w}=\stab^{\fC_+, T\calB, \na_-}_w, \quad \stabm_w=\stab^{\fC_-, T^*\calB, \na_+}_w.\] We list some basic properties of stable bases (see  \cite[Proposition 1]{OS16} and \cite[Lemma 2.2]{SZZ17}): 
\begin{enumerate}
\item The duality: $\langle \stab^{\fC, T^{1/2}, \na}_w, \stab^{-\fC, T^{1/2}_{\opp}, -\na}_v\rangle =\de_{v,w}\in K_T(T^*\calB)$. 
\item $\stabp_w|_w=q^{-\ell(w)/2}\prod_{\be>0, w\be<0}(q-e^{w\be})\prod_{\be>0, w\be>0}(1-e^{w\be})$. 
\item $\stabm_w|_w=q^{\ell(w)/2}\prod_{\be>0, w\be<0}(1-e^{-w\be})\prod_{\be>0, w\be>0}(1-qe^{-w\be})$. 
\end{enumerate}

\subsection{The action of the affine Hecke algebra} Let $\bbH$ denote the affine Hecke algebra of $G$. As a vector space, it is $H(W)\otimes \bbZ[\La]$, where $H(W)$ is the finite Hecke algebra, whose quadratic relation is $(T_{s_\alpha}+1)(T_{s_\alpha}-q)=0$ for a simple root $\alpha$. Recall the famous theorem of Kazhdan--Lusztig \cite{KL} and Ginzburg \cite{CG} 
\begin{equation}\label{equ:KLG}
\bbH\simeq K_{G\times \bbC^*}(Z),
\end{equation}
where $Z=T^*\calB\times_\calN T^*\calB$ is the Steinberg variety, and the right hand side is endowed with the convolution algebra structure (see \cite[Chapter 5]{CG}). The isomorphism is constructed as follows. We use the notations in Section \ref{sec:graded}. Let $\pi_i:T^*_{Y_\al}\to T^*\calB$ are the two projections. Let $\calO_\triangle$ be the structure sheaf of the diagonal $\triangle(T^*\calB)\subset T^*\calB\times T^*\calB$. Then the isomorphism in  \eqref{equ:KLG} sends the simple generator $T_{s_\alpha}\in H(W)$ to $-[\calO_\triangle]-[\calO_{T^*_{Y_\al}}\otimes \pi_2^*\calL_\alpha]\in K_{G\times \bbC^*}(Z)$ and $e^\lambda\in \bbZ[\La]$ to $\calO_\triangle(\lambda)$; see \cite{R08}.

By convolution, $K_{G\times \bbC^*}(Z)$ acts on $K_T(T^*\calB)$; see \cite[Chapter 5]{CG}. Since the kernel defining $T_\alpha$ is not symmetric. The left and right convolution actions are different. Following \cite{SZZ17}, we use $T_\al$ (resp. $T_\al'$) to denote the left (resp. right) convolution action of $T_\alpha$. From \cite[Lemma 3.4]{SZZ17}, these two operators are adjoint to each other:
\[
\langle T_\al(\calF), \calG\rangle =\langle \calF, T_\al'(\calG)\rangle, \quad \forall \calF, \calG\in K_T(T^*\calB). 
\]

One of the main results in \cite{SZZ17} is the computation of the affine Hecke algebra action on the stable bases. More precisely, we have
\begin{thm}\cite[Proposition 3.3, Theorem 3.5]{SZZ17} \label{thm:actiononstable}
	Let $\alpha$ be a simple root. Then
	\begin{align}
	\label{eq:Taction}T_\al(\stabm_{w})&=\left\{ \begin{array}{ll}(q-1)\stabm_{w}+q^{1/2}\stabm_{ws_\al}, & \text{ if }ws_\al<w;\\
	q^{1/2}\stabm_{ws_\al}, &\text{ if }ws_\al>w.\end{array}\right.\\
	T'_\al(\stabp_{w})&=\left\{\begin{array}{ll}(q-1)\stabp_{w}+q^{1/2}\stabp_{ws_\al}, &\text{ if }ws_\al<w;\\
	q^{1/2}\stabp_{ws_\al}, &\text{ if }ws_\al>w.\end{array}\right.
	\end{align}
	In particular,
	\[\stabm_{w}=q^{\frac{\ell(w_0w)}{2}}T^{-1}_{w_0w}(\stabm_{w_0}), \textit{\quad and \quad} \stabp_{w}=q^{-\frac{\ell(w)}{2}}T'_{w^{-1}}(\stabp_{id}).\]
\end{thm}
In the proof of this theorem, an elementary but essential method called rigidity was used. Its simplest form says the following. If $p(z)\in \bbC[z^{\pm 1}]$, then 
\[
p(z) \text{ is bounded as }z\to \pm \infty \quad \Longleftrightarrow \quad  p(z) \text{ is a constant}. 
\]
More applications of this method can be found in the survey by Okounkov \cite{O}. As a immediate corollary, we obtain recursive formulas for  localizations of  stable bases \cite[Proposition 3.6]{SZZ17}. This plays an important role in the identification between the stable basis and motivic Chern classes of  Schubert cells (see Section \ref{sec:motivicChern}), which was used to prove some conjectures of Bump, Nakasuji and Naruse in representations of $p$-adic Langlands dual groups; see \cite{AMSS19}.

\section{Root polynomials  and restriction formulas}
Since localization is such a powerful tool in calculations, it is important to study the localization of the stable basis. The cohomology case was done in \cite{Su17}. For the K-theoretic case, we use the root polynomial method, which works better in an algebraic construction of the K-theory stable basis. 
\subsection{Algebraic construction of stable bases}
Firstly, we can realize stable basis in terms  of Kostant and Kumar's twisted group algebra method. Recall that $S=\bbZ[q^{1/2}, q^{-1/2}][\La]$. Denote $x_\alpha:=1-e^{-\al}$.  Let  $Q=\Frac(S)$, $Q_W=Q\rtimes \bbZ[W]$ with a left $Q$-basis $\de_w, w\in W$ with product formula 
\[
p\de_w \cdot p'\de_{w'}=pw(p')\de_{ww'}, ~w,w'\in W, p,p'\in Q. 
\]
For each simple root $\al_i$, define the Demazure-Lusztig elements in $Q_W$ as follows: 
\[
\taup_i=\frac{q-1}{1-e^{\al_i}}+\frac{q-e^{\al_i}}{1-e^{-\al_i}}\de_{s_i}, \quad \taum_i=\frac{q-1}{1-e^{\al_i}}+\frac{1-qe^{-\al_i}}{1-e^{\al_i}}\de_{s_i}.
\]
Both of the two elements satisfy the usual defining relations for the Demazure-Lusztig operators, so together with $S$, $\taup_i$ (resp. $\taum_i$) for $i=1,...,n$ generate a subalgebra of $Q_W$ that is isomorphic to the  affine Hecke algebra $\bbH$. Moreover, $\taupm_w, w\in W$ form bases of $Q_W$, and we have
\begin{equation}
\de_w=\sum_{v\le w}b^\pm_{w,v}\taupm_{v}, \quad b^\pm_{w,v}\in Q.
\label{eq:tautode}\end{equation}
We will show that the coefficients $b^{\pm}_{w,v}$ can be computed by root polynomials below, and they are closely related to the restriction formula of stable bases.

Define the dual of the left $Q$-module $Q_W$, $Q_W^*:=\Hom_Q(Q_W, Q)\cong \Hom(W, Q)$. It has a standard basis $f_w, w\in W$. This is dual to the basis $\de_w$. The module $Q_W^*$ is a commutative ring via the product $f_wf_v=\de_{w,v}f_w$, and the identity is denoted by $\unit:=\sum_{v\in W}f_v$. Geometrically, $Q_W^*$ is isomorphic to  $Q\otimes_SK_T(T^*\calB)$ as commutative rings, and the basis $f_w$ goes to the normalized fixed point basis $\frac{\iota_{w*}1}{\prod_{\alpha>0}(1-e^{w\alpha})(1-qe^{-w\alpha})}$.

The ring $Q_W$ acts on $Q_W^*$ as follows:
\[
(z\bullet f)(z')=f(z'z), \quad z,z'\in Q_W, ~f\in Q_W^*.
\]
The action is $Q$-linear (hence also $S$-linear). 

We denote some elements in $S$:
\begin{gather*}
x_{w_0}=\prod_{\al>0}(1-e^{-\al}), \quad x_{-w_0}=w_0(x_{w_0}).
\end{gather*} 
By comparing the $T_{\al_i}, T'_{\al_i}$ action with $\taupm_{\al_i}$, and using the normalization condition of stable bases,  we obtain that
\begin{thm}\cite[Theorem 5.4]{SZZ17} Via the isomorphism $Q\otimes_S K_T(T^*\calB)\cong Q_W^*$, $T_{\al_i}$ and $T_{\al_i}'$ coincide with $\taum_i, \taup_i$, respectively, and we have
\begin{align*}
\stabp_w=q^{-\ell(w)/2}\taup_{w^{-1}}\bullet (x_{-w_0}f_e), \quad \stabm_w=q^{\ell(w_0)}q^{-\ell(w)/2}(\taum_{w_0w})^{-1}\bullet (x_{-w_0}f_{w_0}).
\end{align*}
\end{thm}

\subsection{Root polynomials and restriction formula}
Let $Q^x=Q$ when variables of $Q$ are denoted by $x_\la=1-e^{-\la}$, and variables of $Q^y=Q$ will be $y_\la=1-e^{-\la}$. Let $\hat Q=Q^y\otimes_{R} Q^x$ and $\hat Q_W=Q^y\otimes Q^x_W $ so that elements of $Q^y$ commute with elements of $Q^x_W=Q^x\rtimes \bbZ[W]$. The free $\hat Q$-module $\hat Q_W$ has basis $\{\de_w^x\}_{w\in W}$. There is a ring homomorphism of $\ev:\hat Q\mapsto Q^x, y_\la\otimes x_\mu\mapsto x_\la x_\mu$, which induces a left $\hat Q$-module structure on $Q^x_W$. Moreover, it induces a left $\hat Q$-module homomorphism $\ev: \hatQ=Q^y\otimes_RQ_W^x\mapsto Q_W^x$, satisfying 
\[
\ev(y\hat z z)=\ev(y)\ev(\hat z)\ev(z), \quad y\in Q^y, \hat z\in \hat Q_W, z\in Q^x_W.
\]

We now define the root polynomials. For each $w=s_{i_1}\cdots s_{i_l}$, denote $\be_j=s_{i_1}\cdots s_{i_{j-1}}\al_j$, and define 
\[
R_w^\pm=\prod_{j=1}^lh_{i_j}(\be_j)\in \hat Q_W, \text{ where } h_i(\beta)=\taupm^x_i-\frac{q-1}{y_{-\be}}\in \hatQ_W. 
\]
Express $R^\pm_w=\sum_{v\le w}K_{v,w}^{\taupm}\taupm^x_v$. From the braid relations of $\taupm_i^x$, we know that $K^{\taum}_{v,w}=K^{\taup}_{v, w}\in Q^y$.
Moreover, we have
\begin{thm}\cite[Theorem 6.3, 6.5]{SZZ17}
For any $w\in W$, we have
\begin{enumerate}
\item $\ev(R^{\taup}_w)=(\underset{\al>0, w^{-1}\al<0}\prod\frac{q-e^\al}{1-e^{-\al}})\de_w^x$, $\ev(R^{\taup}_w)=(\underset{\al>0, w^{-1}\al<0}\prod\frac{1-qe^{-\al}}{1-e^\al})\de_w^x$.
\item $K^{\taup}_{v,w}=K^{\taum}_{v,w}=(\underset{\al>0, w^{-1}\al<0}\prod\frac{q-e^\al}{1-e^{-\al}})b^+_{w,v}=(\underset{\al>0, w^{-1}\al<0}\prod\frac{1-qe^{-\al}}{1-e^\al})b^-_{w,v}$, where $b^\pm_{w,v}$ were defined in \eqref{eq:tautode}. 
\item $\stabm_w=q^{-\ell(w_0w)}\sum_{v\ge w}\underset{\al>0, v^{-1}\al<0}\prod(1-e^\al)\underset{\al>0, v^{-1}\al>0}\prod(1-qe^{-\al})K^\taupm_{w,v}f_v$. 
\end{enumerate}
\end{thm}
Restricting both sides of the third equation to the fixed point $vB$, we get the formula for $\stabm_w|_v$.

\begin{exmp}We consider the case of $SL_3$, in which case there are two simple roots $\al_1,\al_2$ and $W=S_3$.  Let $w=s_1s_2$ with $\be_1=\al_1, \be_2=\al_1+\al_2$. Then 
\begin{gather*}
R^\taupm_w=h_1(\al_1)h_{2}(\al_1+\al_2)=(\taupm^x_{s_1}-\frac{q-1}{y_{-\al_1}})(\taupm^x_{s_2}-\frac{q-1}{y_{-\al_1-\al_2}})\\
=\taupm^x_{s_1s_2}-\frac{q-1}{y_{-\al_1}}\taupm^x_{s_2}-\frac{q-1}{y_{-\al_1-\al_2}}\taupm^x_{s_1}+\frac{(q-1)^2}{y_{-\al_1}y_{-\al_1-\al_2}}. 
\end{gather*}
So 
\[
\ev(R_w^\pm)=\taupm^x_{s_1s_2}-\frac{q-1}{x_{-\al_1}}\taupm^x_{s_2}-\frac{q-1}{x_{-\al_1-\al_2}}\taupm^x_{s_1}+\frac{(q-1)^2}{x_{-\al_1}x_{-\al_1-\al_2}}. 
\]
For instance, expressing $\stabm_{s_1}$ as a linear combination of $f_v, v\in S_3$, then the coefficient at $s_1s_2$ is equal to 
\begin{align*}&q^{-\ell(w_0s_1)}\underset{\al>0, s_2s_1\al<0}\prod(1-e^\al)\underset{\al>0, s_2s_1\al>0}\prod(1-qe^{-\al})K^\taupm_{s_1, s_1s_2}
\\=&q^{-3}(1-e^{\al_1})(1-e^{\al_1+\al_2})(1-qe^{-\al_2})(-\frac{q-1}{1-e^{\al_1+\al_2}})\\
=&-q^{-3}(q-1)(1-e^{\al_1})(1-qe^{-\al_2}). \end{align*}
\end{exmp}

\section{Motivic Chern classes of the Schubert cells}\label{sec:motivicChern}
In this section, we recall briefly the relation between the stable basis and motivic Chern classes for the Schubert cells, see \cite{AMSS19, FRW}.

\subsection{Definition of motivic Chern classes}
Let us first recall the definition of the motivic Chern classes, following \cite{BSY, AMSS19, FRW}. Let $X$ be a quasi-projective complex $A$-variety. Let ${G}_0^A(var/X)$ be the free abelian group generated by symbols $[f: Z \to X]$ where $Z$ is a quasi-projective $A$-variety and $f: Z \to X$ is a $A$-equivariant morphism modulo the usual additivity relations 
\[[f: Z \to X] = [f: U \to X] + [f:Z \setminus U \to X]\] for $U \subset Z$ an open $A$-invariant subvariety. Then there exists a unique natural transformation 
\[MC_y: G_0^A(var/X) \to K_A(X)[y]\]
such that it is functorial with respect to $A$-equivariant proper morphisms of non-singular, quasi-projective varieties, and if $X$ is smooth,
\[ MC_y([id_X: X \to X]) = \lambda_y(T^*X) = \sum y^i [\wedge^i T^*X] \in K_A(X)[y] \/. \] Here $y$ is a formal variable. The existence in the non-equivariant case was proved in \cite{BSY}. The equivariant case is established in \cite{AMSS19, FRW}.

\subsection{Motivic Chern classes of the Schubert cells}
Let $X=G/B$ be the complete flag variety. Recall $X(w)^\circ=BwB/B\subset X$ and $Y(w)^\circ=B^-wB/B$ are the Schubert cells in $X$, where $B^-$ is the opposite Borel group. Recall the Serre duality functor on $X$
\[\calD(\calF)=R\calH om(\calF, \omega^\bullet_{X}),\]
where $\calF\in K_A(X)$ and $\omega^\bullet_{X}=\calL_{2\rho}[\dim X]$ is the dualizing complex of $X$. We extend it to $K_A(X)[y,y^{-1}]$ by sending $y$ to $y^{-1}$. Let $i:X\hookrightarrow T^*\calB$ be the inclusion. Then the relation between the motivic Chern classes and the stable basis is
\begin{thm}\cite{AMSS19}\label{thm:mot} For any $w\in W$, we have
\[q^{-\frac{\ell(w)}{2}}\calD(i^*\stab_+(w))=MC_{-q^{-1}}([X(w)^\circ\hookrightarrow X])\in K_A(X)[q,q^{-1}],\]
and 
\[q^{\frac{\ell(w)}{2}-\dim G/B}i^*(\stab_-(w))\otimes [ \omega^\bullet_{X}] 
=MC_{-q^{-1}}([Y(w)^\circ\hookrightarrow X])\in K_A(X)[q,q^{-1}].\]
\end{thm}

It will be proved in \cite{AMSSa} that 
\[ (-q)^{- \dim G/B} i^*( {gr(i_{w!}\mathbb{Q}_{Y(w)^\circ}^H)) \otimes [\omega^\bullet_X]} = MC_{-q^{-1}}([Y(w)^\circ\hookrightarrow X]) \/, \]
where $i_w: Y(w)^\circ\hookrightarrow X$ is the inclusion, and {$gr(i_{w!} \mathbb{Q}_{Y(w)^\circ}^H)$ is the associated graded (or $\bbC^*$-equivariant)} sheaf on the cotangent bundle of $X$ determined by the shifted mixed Hodge module $\mathbb{Q}_{Y(w)^\circ}^H$; see \cite{T87, AMSSa}. Since $i^*$ is an isomorphism, we get
\[q^{\frac{\ell(w)}{2}}\stab_-(w)=(-1)^{\dim G/B}gr(i_{w!}\mathbb{Q}_{Y(w)^\circ}^H).\]
This is the K-theoretic analogue of Theorem \ref{thm:char}.

Theorem \ref{thm:mot} is used in \cite{AMSS19} to prove some conjectures of Bump, Nakasuji and Naruse (\cite{BN17, N14}) about unramified principal series of the $p$-adic Langlands dual group. In the proof, we need to build a relation between the stable basis (or motivic Chern classes of Schubert cells) and unramified principal series. This is explained in the next section.

\section{Unramified principle series of $p$-adic Langlands dual group}
In this section, we relate stable bases with certain basis in the Iwahori-invariants in the unramified principal series representation of the Langland dual group. 

\subsection{The two bases in Iwahori-invariants of an unramified principal series representation}
Let $F$ be a finite extension of $\bbQ_p$ \footnote{The result also holds for $F=F_q((t))$.}. Let $\calO_F$ be the ring of integers, with a uniformizer $\varpi\in \calO_F$, and residue field $\bbF_q$. Let $\check{G}$ be the split reductive group over $F$, which is the Langlands dual group of $G$. Let $\check{A}\subset \check{B}$ be the corresponding dual maximal torus and Borel subgorup. Let $I$ be the Iwahori subgroup. 

Let $\tau$ be an unramified character of $\check{A}$. Recall the principal series representation  $I(\tau):=\Ind^{\check{G}(F)}_{\check{B}(F)}(\tau)$ is the induced representation, which consists of locally constant functions $f$ on $\check{G}(F)$ such that $f(bg)=\tau(b)\de^{1/2}(b)f(g), b\in \check{B}(F)$, where $\de(b)=\prod_{\al>0}|\al^\vee(b)|_F$ is the modulus function on the Borel subgroup.

It is well known that the affine Hecke algebra in  \eqref{equ:KLG} can also be realized as $\bbH=\bbC_c[I\backslash \check{G}(F)/I]$, which is the so-called Iwahori--Hecke algebra. The algebra structure on the latter space is induced by convolution. Again by convolution, the Iwahori--Hecke algebra $\bbH$ acts from the right on the Iwahori-invariants $I(\tau)^I$. 

If the unramified character $\tau$ is regular, then the space $\Hom_{\check{G}(F)}(I(\tau), I(w^{-1}\tau))$ is one-dimensional. It is spanned by an operator $\calA^\tau_w$, called intertwiner. And it is defined as follows
\[
\calA^\tau_w(f)(g)=\int_{\check{N}_w}f(wng)dn,
\]
where $\check{N}_w=\check{N}(F)\cap w^{-1}\check{N}^-(F)w$ with $\check{N}$ (resp. $\check{N}^-$) being the unipotent radical of the Borel subgroup $\check{B}$ (resp. $\check{B}^-$).

There are two bases of $I(\tau)^I$. The first basis $\{\varphi^\tau_w| w\in W\}$, which is called the standard basis, is induced by the following decomposition 
\[\check{G}(F)=\sqcup_{w\in W}\check{B}(F)wI.\] 
I.e., $\varphi^\tau_w$ is characterized by the following conditions:
\begin{enumerate}
\item $\varphi^\tau_w$ is supported on $\check{B}(F)wI$;
\item $\varphi^\tau_w(bwg)=\tau(b)\delta^{1/2}(b)$ for any $b\in \check{B}(F), g\in I$. 
\end{enumerate}
The other basis $\{f_w^\tau| w\in W\}$, which is the so-called Casselman's basis \cite{C80}, is characterized by
\[\calA_v^\tau(f^\tau_w)(1)=\de_{v,w}. \]
These bases played an important role in the computation of the Macdonald spherical function \cite{C80} and the Casselman--Shalika formula for the spherical Whittaker function \cite{CS}.

\subsection{The comparison}
Since $\tau$ is an unramified characte of $\check{A}$, it is equivalent to a point $\tau\in A$. Therefore, we can evaluate the base ring $K_A(pt)$ at $\tau\in A$. We let $\bbC_\tau$ denote this evaluation representation of $K_T(pt)=K_{A\times \bbC^*}(pt)$. For any $\calF\in K_{A\times \bbC^*}(T^*\calB)$, let $\calF_{-\rho}$ denote $\calF\otimes \calL(-\rho)$. 

Now we can state the relation between these bases. The following theorem is a shadow of the two geometric realizations of the affine Hecke algebra 
\[K_{G\times \bbC^*}(Z)\simeq\bbH\simeq \bbC_c[I\backslash \check{G}(F)/I].\]
We refer the readers to \cite[Section 8.2 and 8.3]{SZZ17} for a more precise statement of the following theorem.

\begin{thm}\cite[Theorem 8.4]{SZZ17}\label{thm:kpadic}
There is a unique isomorphism of right $\bbH$-modules 
\[
\Psi: K_{A\times \bbC^*}(T^*\calB)\otimes_{K_{A\times \bbC^*}(\pt)}\bbC_\tau\to I(\tau)^I,
\]
such that the equivariant parameter $q$ is mapped to the cardinality of the residue field $\bbF_q$ and
\begin{align*}
\frac{q^{\ell(w)}}{\prod_{\be>0, w\be>0}(1-e^{w\be})\prod_{\be>0, w\be<0}(q-e^{w\be})}(\iota_{w*}1)_{-\rho}&\mapsto& f^\tau_w, \\
(\stabm_w)_{-\rho} &\mapsto& q^{-\ell(w)/2}\varphi^\tau_w. 
\end{align*}
\end{thm}
As applications of this theorem, in \cite{SZZ17}, we showed that the left Weyl group action on the K-theory side corresponds to the suitably normalized intertwiners on the representation theory side (cf. Corollary 8.6), and  provided some K-theoretic interpretations of Macdonald's formula \cite{C80} (cf. Theorem 8.8) and the Casselman-Shalika formula \cite{CS} (cf. Theorem 8.11). 

It is an interesting question to study the transition matrix between the two basis in $I(\tau)^I$. Define the following transition matrix coefficients $m_{u,w}$ by
\[\sum_{w\geq u}\varphi_w=\sum_{w\in W} m_{u,w}f_w.\] For instance, the following special case
\[m_{id,w} = \prod_{\alpha>0, w^{-1}\alpha<0} \frac{1 - q^{-1} e^{\alpha}(\tau)}{1 - e^{\alpha}(\tau)} \] 
is the Gindikin-Karpelevich formula. 

The authors of \cite{AMSS19} used Theorem \ref{thm:kpadic}, Theorem \ref{thm:mot} and some functorial properties for the motivic Chern classes to prove the following conjectures of Bump, Nakasuji and Naruse (\cite{BN17, N14}).
\begin{thm}\cite{AMSS19}
For any $u\leq w\in W$, 
\begin{enumerate}
\item
\begin{equation*}
m_{u,w}=\prod_{\alpha>0, u\leq s_\alpha w < w}\frac{1-q^{-1}e^\alpha(\tau)}{1-e^\alpha(\tau)},
\end{equation*}
if and only if the Schubert variety $Y(u)$ is smooth at the torus fixed point $e_{w}$.
\item 
As a function of $\tau\in A$, the product $\prod_{\alpha>0, u\leq s_\alpha w < w}(1-e^\alpha)m_{u,w}$ is analytic on the maximal torus $A$.
\end{enumerate}
\end{thm} 

\section{Wall crossings for the stable bases}
Note that the definition of the stable bases depend on the alcoves in $\fh^*_\bbR$. Stable bases for different alcoves are related by the so-called wall crossing R-matrices; see \cite{OS16}. In this section, we study the wall crossing R-matrices for the Springer resolution, which will be used for the categorification in the next section. The main reference for these two sections is \cite{SZZ19}.

\subsection{Wall crossing matrix}
By  uniqueness of the stable basis, it is immediate to see that for any $\mu\in \La$,
\begin{equation}\label{equ:shift}
\stab^{\fC,T^{1/2},\na+\mu}_y=e^{-y\mu}\calL_\mu\otimes \stab^{\mathfrak{C},T^{1/2},\na}_y.
\end{equation}
Because of this property, instead of crossing all the walls on $H_{\alpha^\vee, n}$, it is enough to just cross the walls on the $0$ hyperplanes $H_{\alpha^\vee, 0}$. From now on, let $\na_1$, $\na_2$ are two alcoves sharing a wall on $H_{\al^\vee, 0}$, and $(\la_1, \al^\vee)>0$ for any $\la_1\in \na_1$.

Another useful fact is (see \cite[Theorem 1]{OS16}):
\[\stab^{\fC, T^{1/2}, \na_2}_y=\left\{\begin{array}{ll} \stab^{\fC, T^{1/2}, \na_1}_y+f^{\na_2\leftarrow \na_1}_y\stab^{\fC, T^{1/2}, \na_1}_{ys_\al}, & \text{ if } ys_\al \prec_\fC y;\\
\stab^{\fC, T^{1/2}, \na_1}_y, & \text{ if }ys_\al \succ_\fC y, \end{array}\right.\]
where $f^{\na_2\leftarrow \na_1}_y\in K_T(\pt)$. 

Then we have
\begin{thm}\cite{SZZ19}\label{thm:crosswall}
For any $y\in W$, we have
\begin{gather*}
\stab_y^{\fC, T\calB, \na_2}=\left\{\begin{array}{ll}\stab^{\fC, T\calB, \na_1}_y+(q^{1/2}-q^{-1/2})\stab^{\fC, T\calB, \na_1}_y, & \text{ if } ys_\al\prec_\fC y;\\
\stab_y^{\fC, T\calB, \na_1}, & \text{ if }ys_\al \succ_\fC y.  \end{array}\right.\\
\stab_y^{-\fC, T^*\calB, \na_2}=\left\{\begin{array}{ll}\stab^{-\fC, T^*\calB, \na_1}_y+(q^{1/2}-q^{-1/2})\stab^{-\fC, T^*\calB, \na_1}_y, & \text{ if } ys_\al\succ_\fC y;\\
\stab_y^{-\fC, T^*\calB, \na_1}, & \text{ if }ys_\al \prec_\fC y.  \end{array}\right.\\
\end{gather*}
\end{thm}
This theorem is firstly proved when $\al$ is a simple root, by computing the action of $T_\al$ and $T_\al'$ on the stable basis $\stab_y^{\fC, T\calB, \na}$ and $\stab^{-\fC, T^*\calB, \na}_y$ for any alcove $\na$. I.e., we use rigidity to obtain a general version of Theorem \ref{thm:actiononstable}. For the non-simple root case, we use the following equality (\cite{AMSS19}) to reduce to the simple root case 
\begin{equation}\label{equ:weyl}
w(\stab_y^{\fC, T\calB, \na})=\stab_{wy}^{w\fC, T\calB, \na},
\end{equation}
where $w$ denotes the left Weyl group action on $K_T(T^*\calB)$ by $w$.

\subsection{Wall crossing and affine Hecke algebra actions}
Combining Theorem \ref{thm:actiononstable} and Theorem \ref{thm:crosswall}, we get the following general formulae, which show that the wall crossing matrices and the affine Hecke algebra action (see Section \ref{sec:hecke}) are compatible. 
\begin{thm}\cite{SZZ19}\label{thm:wallcrossingHecke}
	For any $x, y\in W$, we have 
	\begin{eqnarray}\label{eq:Neg}
    \stab^{\fC_-, T^*\calB, x\na_+}_y=&q_x^{-1/2}T_x(\stab^{\fC_-, T^*\calB, \na_+}_{yx});\\ 
	\stab^{\fC_+, T\calB, x\na_-}_y=&q_x^{1/2}(T'_{x^{-1}})^{-1}(\stab^{\fC_+, T\calB, \na_-}_{yx}). \label{eq:Pos}
	\end{eqnarray}
\end{thm}
Therefore, Theorem \ref{thm:actiononstable}, Equations \eqref{equ:shift} and \eqref{eq:Pos} determine all the stable basis $\stab^{\fC_+, T\calB, \na}_y$ for the dominant Weyl chamber $\fC_+$. The stable basis for the other chambers can be computed by Equation \eqref{equ:weyl}. Thus all the stable basis element $\stab^{\fC, T\calB, \na}_y$ can be calculated.

For general symplectic resolutions, Bezrukavnikov and Okounkov (\cite{O15, BO}) conjecture that the representation coming from the derived equivalence is isomorphic to the monodromy representation coming from quantum cohomology. The monodromy matrices of the quantum connection of $T^*\calB$ is computed in \cite{BMO, C}. The above theorem shows that the monodromy matrices coincide with the wall R-matrices for the K stable basis. The relation to derived equivalences is explained in the next section.

\section{Categorification and localization in positive characteristic}
In this final section, we give a categorification of the stable bases, and study its relation with representation of $\fg$ over positive characteristic fields under the localization equivalence of Bezrukavnikov, Mirkovi\'{c} and Rumynin \cite{BMR1, BMR2}.

\subsection{Categorification of the stable basis via affine Braid group action}Let $B_{\aff}$ (resp. $B'_{\aff}$) be the affine braid group (resp. the extended affine braid group) with generators $\wt s_\al, \al\in I_{\aff}$. 

Let $G_\bbZ$ be a split $\bbZ$-form of the complex algebraic group $G$ and $A_\bbZ\subset B_\bbZ\subset G_\bbZ$ be the maximal torus and a Borel subgroup, respectively. Let $T_\bbZ=A_\bbZ\times_\bbZ(\bbG_m)_\bbZ$. Bezrukavnikov and Riche constructed an extended affine braid group action on $D^{\ob}_{T_\bbZ}(T^*\calB_\bbZ)$ (see \cite{BR, R08}), denoted by $J^{\oR}_{\wt w}, \wt w\in B'_{\aff}$ (here $\oR$ denotes the right action). Inspired by Theorem \ref{thm:wallcrossingHecke}, we give the following definition.

\begin{defn} Let $\la\in \Lambda_\bbQ$ be regular. We define $\fstab^\bbZ_\la(y)\in D^{\ob}_{T_\bbZ}(T^*\calB_\bbZ), y\in W$ as follows:
\begin{align*}
\fstab^{\bbZ}_{\lambda_0}(id)=\calL_{-\rho}\otimes\calO_{T^*_{id}\calB_\bbZ}, &\quad \la_0\in \na_-,  \\
\fstab^{\bbZ}_{\lambda_0}(y)=J_{\wt{y}}^{\oR}\fstab^{\bbZ}_{\lambda_0}(id), & \quad \la_0\in \na_-, \\ 
\fstab^\bbZ_{\lambda}(y)=(J_{\wt{x}}^{\oR})^{-1}\fstab^\bbZ_{\lambda_0}(yx), &\quad y,x\in W, x\lambda_0=\lambda, \la_0\in \na_-,\\ 
\fstab^\bbZ_{\lambda}(y)=e^{-y\mu}J^{\oR}_{\mu}\fstab^\bbZ_{\lambda_1}(y), &\quad y\in W, \mu+\lambda_1=\lambda, \mu\in Q, \lambda_1\in W\nabla_-. 
\end{align*}
\end{defn}
As an immediate corollary of Theorem \ref{thm:actiononstable}, Theorem \ref{thm:wallcrossingHecke}, Equation \eqref{equ:shift} and \cite[Theorem 1.3.1, Proposition 1.4.3, Theorem 1.6.1]{BR}, we have the following theorem, which gives a categorification of the stable basis.
\begin{thm}\cite{SZZ19}Applying the derived tensor $-\otimes^L_\bbZ\bbC$ to $\fstab^\bbZ_\la(w)\in D^{\ob}_{T_\bbZ}(T^*\calB_\bbZ)$, and taking the class in the Grothendieck group, we get $\calL_{-\rho}\otimes \stab_w^{\fC_+, T\calB, \la}\in K_T(T^*\calB)$. 
\end{thm}

\subsection{Verma modules in positive characteristic}
In this section, we briefly sketch the relation between the K-theory stable basis and Verma modules for Lie algebras over positive characteristic fields. We refer the readers to \cite[Section 9]{SZZ19} for the details.

We consider the level-$p$ configuration of $\rho$-shifted affine hyperplanes, that is, $H^p_{\al^\vee, n}=\{\la\in \La_\bbQ |\langle\al^\vee, \la+\rho\rangle=np\}$. Let $A_0$ be the fundamental alcove, i.e., it contains $(\ep-1)\rho$ for small $\ep>0$. Let $W$ act on $\La$ via the level-$p$ dot action,  $w:\la\mapsto w\bullet \la=w(\la+\rho)-\rho$.   

Let $k$ be an algebraically closed field of characteristic $p$, and $p$ is greater than the Coxeter number. For any $k$-variety $X$, let $X^{(1)}$ be the Frobenius twist. Let $U(\fg_k)$ be the universal enveloping algebra of $\fg_k$ with the Frobenius center $\calO(\fg_k^{*(1)})$ and the Harish-Chandra center $\calO(\fh^*_k/(W, \bullet))$. Let $\la\in \fh^*_k$ be regular (i.e., does not lie on any hyperplane $H^p_{\al^\vee, n}$) and integral (i.e., belongs to the image of the derivative $d:\La\to \fh^*_k$). Let $U(\fg_k)^\la$ be the quotient of $U(\fg_k)$ by the central ideal corresponding to $W\bullet d\la\in \fh^*_k/(W, \bullet)$. Let $\Mod_\chi(U(\fg_k)^\la)$ be the category of finitely generated $U(\fg_k)^\la$-modules on which the Frobenius center $\calO(\fg^{*(1)}_k)$ acts by the generalized character $\chi\in \fg^{*(1)}_k$. 

Let $\calD^\la$ be the ring of $\calL_\la$-twisted differential operators on $\calB$, and $D^{\ob}(\Coh_\chi \calD^\la)$ be the full subcategory of coherent $\calD^\la$-modules that are set-theoretically supported on $\calB^{(1)}_\chi$ as a coherent sheaf on $T^*\calB^{(1)}_k$, where $\calB^{(1)}_\chi$ is the Springer fiber. Then the global section functor $R\Gamma_{\calD^\la, \chi}: D^{\ob}(\Coh_\chi\calD^\la)\to D^{\ob}(\Mod_\chi U(\fg_k)^\la)$ is an equivalence according to \cite[Theorem 3.2]{BMR1}, whose inverse is denoted by $\fL^{\lambda_0}$. Let $T^*\calB^{(1)\wedge}_{\chi}$ be the completion of $T^*\calB^{(1)}$ at $\calB^{(1)}_\chi$. Then for all integral $\lambda\in\fh^*$, the Azumaya algebra $\calD^\lambda$ splits on $T^*\calB^{(1)\wedge}_{\chi}$; see \cite[Theorem~5.1.1]{BMR1}. In particular, there is a Morita equivalence \[\Coh_{\calB_{\chi}^{(1)}}(T^*\calB^{(1)})\cong\Coh_\chi\calD^\lambda.\]
The equivalence above depends  on the choice of a splitting bundle  of $\calD^\lambda$ \cite[Remark~5.2.2]{BMR1}. As has been done in \cite{BR,BM}, for $\lambda_0\in A_0$, we normalize the choice of the splitting bundle by \cite[Remark 1.3.5]{BMR2}. In this section, we denote this splitting bundle by $\calE^s$. That is, $\calE nd_{T^*\calB^{(1)\wedge}_{\chi}}(\calE^s)\cong \calD^{\lambda_0}|_{T^*\calB^{(1)\wedge}_{\chi}}$. This bundle fixes the equivalence \mbox{$\Coh_{\calB_{\chi}^{(1)}}(T^*\calB^{(1)})\cong\Coh_\chi\calD^{\lambda_0}$}. 
Composing this equivalence with $\fL^{\lambda_0}$, we have 
\begin{equation}\label{eqn:gamma_0}
\gamma_\chi^{\lambda_0}: D^{\ob}\Mod_\chi(U(\fg_k)^{\lambda_0})\cong D^{\ob}\Coh_{\calB_\chi^{(1)}}(T^*\calB^{(1)}),
\end{equation}
which is also referred to as the localization equivalence.  

Let $U(\fg_k)^{\lambda_0}_\chi$ be the completion of $U(\fg_k)^{\lambda_0}$ at the central character $\chi\in\calN_k^{(1)}$. Then, $End_{T^*\calB^{(1)\wedge}_{\chi}}(\calE^s)\cong U(\fg_k)^{\lambda_0}_\chi$. The localization functor $\gamma_\chi^{\lambda_0}$ can then be written as $\otimes_{U(\fg_k)^{\lambda_0}_\chi}\calE^s$. We also have the completed version of the equivalence $\gamma^{\lambda_0}_\chi:  D^{\ob}\Mod(U(\fg_k)^{\lambda_0}_\chi)\cong D^{\ob}\Coh(T^*\calB^{(1)\wedge}_{\chi})$ \cite[Remark~2.5.5]{BR}.  

From now on we consider $\chi=0$. In this case, we abbreviate $\gamma_\chi^{\lambda_0}$ simply as $\gamma^{\lambda_0}$, and the twisted Springer fiber $\calB_{\chi}^{(1)}$ is the zero section $\calB^{(1)}\subseteq T^*\calB^{(1)}$.

The bundle $\calE^s$ on $T^*\calB^{(1)\wedge}_{0}$ has a natural $T_k$-equivariant structure \cite[\S~5.2.4]{BM}, where $T_k=A_k\times k^*. $
Taking the ring of endomorphisms, we get a $T_k$-action on $U(\fg_k)_0^{\lambda_0}$ compatible with that on $\calN^{(1)}$. In particular, the 
$\Gm\subseteq T_k$-action provides a non-negative grading on $U(\fg_k)_0^{\lambda_0}$,  referred to as the Koszul grading.
The localization equivalences above can be made into equivalences of equivariant categories.
Let $\Mod_0^{\gr}(U(\fg_k)^{\lambda_0},A_k)$ be the category of finite-dimensional Koszul-graded modules of $U(\fg_k)^{\lambda_0}_0$, which are endowed with compatible actions of the subgroup $A_k\subseteq G_k$. The compatibility is in the sense that  the action of $A_k$ differentiates to the action of the subalgebra $\fh_k\subseteq \fg_k$. Then, we have  \cite[Theorem~1.6.7]{BM}
\[\gamma^{\lambda_0}: D^{\ob}\Mod_0^{\gr}(U(\fg_k)^{\lambda_0},A_k)\cong D^{\ob}_{T_k}\Coh_{\calB_k^{(1)}}(T^*\calB_k^{(1)}),\]
 \[\gamma^{\lambda_0}: D^{\ob}\Mod^{\gr}(U(\fg_k)_0^{\lambda_0},A_k)\cong D^{\ob}_{T_k}\Coh(T^*\calB^{(1)\wedge}_{0}).\]
By the functor of taking finite vectors, we get (\cite[Theorem~5.1.1]{BM})
\[\gamma^{\lambda_0}: D^{\ob}\Mod^{\gr}(\calA^{\lambda_0},A_k)\cong D^{\ob}_{T_k}\Coh(T^*\calB^{(1)}_k).\]
Here $\calA^{\lambda_0}$ is an $\calO_{\calN^{(1)}}$-algebra endowed with a compatible $A_k\times(\Gm)_k$-action, and it has the property that $(\calA^{\lambda_0})^\wedge_0\cong U(\fg_k)_0^{\lambda_0}$. Using the correspondence given by taking completion and taking finite vectors, we will freely pass between $\Mod^{\gr}(U(\fg_k)_0^{\lambda_0},A_k)$ and $\Mod^{\gr}(\calA^{\lambda_0},A_k)$; similarly for $D^{\ob}_{T_k}\Coh(T^*\calB^{(1)\wedge}_{0})$ and $D^{\ob}_{T_k}\Coh(T^*\calB^{(1)})$.

For any $\lambda$ in the $W_{\aff}'$-orbit of $\lambda_0$, we can define the localization functor 
\[\gamma^{\lambda}: D^{\ob}\Mod^{\gr}(\calA^{\lambda_0},A_k)\cong D^{\ob}_{T_k}\Coh(T^*\calB^{(1)}_k).\]
by precomposing with the affine braid group action functors; see \cite[Section 9.3]{SZZ19}.

For the Lie algebra $\fb_k$ and $\lambda\in\Lambda$, recall the Verma module $Z^{\fb}(\lambda):=U(\fg)\otimes_{U(\fb)}k_\lambda$. Then we have
\begin{thm}\cite{SZZ19}
Let $k$ be an algebraically closed field of characteristic $p$, and $p$ is greater than the Coxeter number. Assume $\lambda$ to be regular and integral, then  in $D^{\ob}_{T_k}(T^*\calB_k^{(1)})$, we have isomorphisms	\[e^{\rho}\fstab^k_{-\frac{\lambda+\rho}{p}}(y)\cong \gamma^{\lambda} Z^{\fb}(y\bullet \lambda+2\rho),\]
where $\fstab^k_{-\frac{\lambda+\rho}{p}}(y)=\fstab^{\bbZ}_{-\frac{\lambda+\rho}{p}}(y)\otimes^L_{\bbZ}k$.
\end{thm}
To prove this theorem, one needs to use the affine braid group action on $\Mod_\chi U(\fg_k)^\la$ constructed in \cite{BMR2}, which iterates the Verma modules \cite{Janz},  and also its compatibility with the localization equivalences $\gamma^\la$. Together with the iterative definition of $\fstab^k_\bbZ(y)$, the theorem will follow.



\newcommand{\arxiv}[1]
{\texttt{\href{http://arxiv.org/abs/#1}{arXiv:#1}}}
\newcommand{\doi}[1]
{\texttt{\href{http://dx.doi.org/#1}{doi:#1}}}
\renewcommand{\MR}[1]
{\href{http://www.ams.org/mathscinet-getitem?mr=#1}{MR#1}}


\begin{thebibliography}{BK}

\bibitem[AO16]{AO} M. Aganagic and A. Okounkov, {\em Elliptic stable envelope}, arXiv preprint arXiv:1604.00423

\bibitem[AM16]{AM16} P. Aluffi and L. Mihalcea, {\em Chern--Schwartz--MacPherson classes for Schubert cells in flag manifolds}, Compositio Mathematica, Vol 152, No. 12, 2603--2625, 2016

\bibitem[AMSS17]{AMSS17} P. Aluffi, L. Mihalcea, J. Sch{\"u}rmann and C. Su, {\em Shadows of characteristic cycles, Verma modules, and positivity of Chern-Schwartz-MacPherson classes of Schubert cells}, \arxiv{1709.08697}


\bibitem[AMSS19]{AMSS19} P. Aluffi, L. Mihalcea, J. Sch{\"u}rmann and C. Su, {\em Motivic Chern classes of Schubert cells with applications to Casselman's problem}, preprint \arxiv{1902.10101}. 

\bibitem[AMSS]{AMSSa} P. Aluffi, L. Mihalcea, J. Sch{\"u}rmann and C. Su, {\em Equvivariant Motivic Chern classes via cotangent bundles}, in preparation.



\bibitem[BB81]{BB} A. Beilinson and J. Bernstein, {\em Localisation de g-modules}, C. R. Acad. Sci. Paris S\'er. I Math. 292
(1981), no. 1, 15--18. 



\bibitem[Be16]{B16} R. Bezrukavnikov,
\textit{On two geometric realizations of an affine Hecke algebra},
Publications mathématiques de l'IHéS, 2016, Volume 123, no. 1, Page 1-67.

\bibitem[BM13]{BM} R. Bezrukavnikov and  I. Mirkovi\'{c}, {\em Representations of semisimple Lie algebras in prime characteristic and the noncommutative Springer resolution}. With an appendix by Eric Sommers. Ann. of Math. (2) 178 (2013), no. 3, 835--919. 

\bibitem[BMR06]{BMR2} R. Bezrukavnikov, I. Mirkovi\'{c}, and D.~Rumynin, {\em Singular localization and intertwining functors for reductive lie algebras in prime characteristic}. Nagoya J. Math. \textbf{184} (2006), 1-55. 

\bibitem[BMR08]{BMR1} R. Bezrukavnikov, I. Mirkovi\'{c}, and D.~Rumynin, {\em Localization of modules for a semisimple Lie algebra in prime characteristic,}
With an appendix by Bezrukavnikov and Simon Riche.
Ann. of Math. (2) 167 (2008), no. \textbf{3}, 945-991. 


\bibitem[BO]{BO} R. Bezrukavnikov, and A. Okounkov, {\em in preparation}


\bibitem[BR13]{BR} R. Bezrukavnikov and S. Riche, {\em Affine braid group actions on derived categories of Springer resolutions.} Ann. Sci. \'{E}c. Norm. Sup\'er. (4) \textbf{45} (2012), no. 4, 535--599 (2013). 

\bibitem[B99]{B99} S. Billey, {\em Kostant polynomials and the cohomology ring for G/B}, Duke mathematical journal, Vol 96, No. 1, 205--224, 1999

\bibitem[BSY10]{BSY} J. Brasselet, J. Sch{\"u}rmann and S. Yokura, {\em Hirzebruch classes and motivic Chern classes for singular spaces}, Journal of Topology and Analysis, Vol 2, No. 01, 1--55, 2010

\bibitem[BK81]{BK81} J. Brylinski and M. Kashiwara, {\em Kazhdan-Lusztig conjecture and holonomic systems}, Inventiones mathematicae, Vol 64, No. 3, 387--410, 1981

\bibitem[BMO09]{BMO} A. Braverman, D. Maulik, and A. Okounkov, {\em Quantum cohomology of the Springer resolution} Advances in Mathematics, Vol 227, No 1. 421--458, 2011


\bibitem[BN17]{BN17} D. Bump and M. Nakasuji, {\em Casselman's basis of Iwahori vectors and Kazhdan-Lusztig polynomials}, Canadian Journal of Mathematics, 1--16, 2017

\bibitem[C80]{C80} W. Casselman, {\em The unramified principal series of {p}-adic groups. I. The spherical function}, Compositio Math, Vol 40, 387--406, 1980

\bibitem[CS80]{CS} W. Casselman and J. Shalika, {\em The unramified principal series of $p$-adic groups. II. The Whittaker function}, Compositio Mathematica, Vol 41, No 2, 207--231, 1980

\bibitem[C05]{C} Cherednik, Ivan {\em Double affine Hecke algebras}, Cambridge University Press, Vol 319, 2005

\bibitem[CG97]{CG} N. Chriss, V. Ginzburg, {\em Representation theory and complex geometry} Birkh{\"a}user, Boston-Basel-Berlin, 1997


\bibitem[FRW18]{FRW} L. Feher, R. Rimanyi and A. Weber, {\em Motivic Chern classes and K-theoretic stable envelopes}, \arxiv{1802.01503}

\bibitem[G85]{G} V. Ginzburg, {\em Deligne-Langlands conjecture and representations of affine Hecke algebras}, preprint 1985




\bibitem[J00]{Janz} J.C.~Janzten, {\em Modular representations of reductive Lie algebras}, Journal of Pure and Applied Algebra 152 (2000) 133--185.


\bibitem[Kal09]{Kal09} D. Kaledin, {\em Geometry and topology of symplectic resolutions}, Algebraic Geometry—Seattle 2005, Vol 2, 595--628, 2009

\bibitem[KL87]{KL} D. Kazhdan and G. Lusztig, {\em Proof of the Deligne-Langlands conjecture for Hecke algebras}, Inventiones mathematicae, vol 87, No 1, 153--215, 1987

\bibitem[LZ14]{LZ14} C. Lenart, K. Zainoulline,
{\em Towards generalized cohomology Schubert calculus via formal root polynomials}, Math. Research Letters 24, No. 3, (2017), 839--877.  


\bibitem[L88]{L88} G. Lusztig, {\em Cuspidal local systems and graded Hecke algebras, I}, Publications Math{\'e}matiques de l'Institut des Hautes {\'E}tudes Scientifiques, Vol 67, 145--202, 1998


\bibitem[M74]{M} R. MacPherson, {\em Chern classes for singular algebraic varieties}, Annals of Mathematics (2), 100: 423--432, 1974

\bibitem[MO12]{MO} D. Maulik and A. Okounkov, {\em Quantum groups and quantum cohomology}, \arxiv{1211.1287}


\bibitem[N14]{N14} H. Naruse, {\em Schubert calculus and hook formula}, Slides at 73rd S{\'e}m. Lothar. Combin., Strobl, Austria, 2014


\bibitem[O06]{O06} T. Ohmoto, {\em Equivariant {C}hern classes of singular algebraic varieties with group actions}, Math. Proc. Cambridge Philos. Soc., No. 1, Vol 140, 115--134, 2006

\bibitem[O15a]{O} A. Okounkov, {\em Lectures on K-theoretic computations in enumerative geometry}, Geometry of moduli spaces and representation theory, 251–380, IAS/Park City Math. Ser., 24, Amer. Math. Soc., Providence, RI, 2017. 

\bibitem[O15b]{O15} A. Okounkov, {\em Enumerative geometry and geometric representation theory}, Algebraic Geometry: Salt Lake City 2015 (Part 1), Vol 97, 419-458, 2018

\bibitem[O18]{O18} A. Okounkov, {\em On the crossroads of enumerative geometry and geometric representation theory}, \arxiv{1801.09818}

\bibitem[OS16]{OS16} A. Okounkov and A. Smirnov, {\em Quantum difference equation for Nakajima varieties}, preprint, (2016). \arxiv{1602.09007} 

\bibitem[R08]{R08} S. Riche, 
{\it Geometric braid group action on derived categories of coherent sheaves},
Representation Theory of the American Mathematical Society, Vol. 12, No. 5, (2008), 131--169.

\bibitem[RTV14]{RTV} 
R. Rimanyi, V. Tarasov and A. Varchenko, {\em Trigonometric weight functions as K-theoretic stable envelope maps for the cotangent bundle of a flag variety}, J. Geom. Phys. 94 (2015), 81–119

\bibitem[RV15]{RV15} R. Rimanyi and A. Varchenko, {\em Equivariant Chern-Schwartz-MacPherson classes in partial flag varieties: interpolation and formulae}, \arxiv{1509.09315}

\bibitem[S65a]{S65a} M. Schwartz, {\em Classes caract\'eristiques d\'efinies par une stratification d'une vari\'et\'e analytique complexe. {I}}, C. R. Acad. Sci. Paris, Vol 260, 3262--3264, 1965
	
\bibitem[S65b]{S65b} M. Schwartz, {\em Classes caract\'eristiques d\'efinies par une stratification d'une vari\'et\'e analytique complexe. {II}}, C. R. Acad. Sci. Paris, Vol 260, 3535--3537, 1965
	

\bibitem[Su17]{Su17} C. Su, {\em Restriction formula for stable basis of the Springer resolution}, Selecta Mathematica, Vol 23, No. 1, 497--518, 2017

\bibitem[SZZ17]{SZZ17} C. Su, G. Zhao, and C. Zhong, {\em On the $K$-theory stable bases of the Springer resolution.} 43 pages, to appear in Annales sc. de l'ENS. \arxiv{1708.08013}

\bibitem[SZZ19]{SZZ19} C. Su, G. Zhao, and C. Zhong, {\em Wall-crossings and a categorification of the K-theory stable bases of the Springer resolutions}, \arxiv{1904.03769} 


\bibitem[T87]{T87} T. Tanisaki, {\em Hodge modules, equivariant {$K$}-theory and {H}ecke algebras}, Publ. Res. Inst. Math. Sci., Vol 23, 1987, No. 5, 841--879


\end{thebibliography}
\end{document}